\documentclass[10pt]{amsart}

\usepackage{hyperref}
\usepackage{amssymb, latexsym, amsmath, amsfonts}
\usepackage{amscd}
\usepackage{curves}

\newcommand{\cal}{\mathcal} 

\theoremstyle{plain} 

\newtheorem{theorem}{Theorem}[section] 
\newtheorem{lemma}[theorem]{Lemma}
\newtheorem{proposition}[theorem]{Proposition}
\newtheorem{corollary}[theorem]{Corollary}

\newtheorem{warning}[theorem]{Warning}

\theoremstyle{remark} 

\newtheorem{remark}[theorem]{Remark}

\newcommand{\Lemmaref}[1]{Lem\-ma~\ref{#1}}
\newcommand{\Propref}[1]{Pro\-po\-si\-tion~\ref{#1}} 
\newcommand{\Thmref}[1]{Theo\-rem~\ref{#1}}
\newcommand{\Rkref}[1]{Re\-mark~\ref{#1}}
\newcommand{\Corref}[1]{Co\-rol\-la\-ry~\ref{#1}}
\newcommand{\Figref}[1]{Fi\-gu\-re~\ref{#1}}

\newcommand{\mraise}[2]{\raisebox{#1}{$#2$}}      

\newcommand{\sub}[1]{\raisebox{-2pt}{$\!_{#1}$} }  
\newcommand{\Sub}[1]{\raisebox{-2pt}{$\!_{\,#1}$} }

\newcommand{\adjustnabla}[1]{\sub{#1}} 

\newcommand{\goth}[1]{\mathfrak{#1}}              
\newcommand{\defemph}[1]{{\sffamily\slshape #1}}  

\renewcommand{\phi}{\varphi}

\newcommand{\ie}{{\it i.e.}~}
\newcommand{\st}{such that}
\newcommand{\Iff}{if and only if}

\newcommand{\PD}{Poin\-car\'e-dual}
\newcommand{\LC}{Levi-Civit\`a}
\newcommand{\SW}{Sei\-berg--Witten}

\newcommand{\SD}{self-dual}
\newcommand{\ASD}{anti-self-dual}
\newcommand{\ac}{al\-most-com\-plex}
\newcommand{\riem}{Riemannian}
\newcommand{\herm}{Hermitian}
\newcommand{\kah}{K\"ahler}

\newcommand{\str}{struc\-ture}

\newcommand{\aR}{\mathbb{R}}
\newcommand{\Zi}{\mathbb{Z}}
\newcommand{\Ce}{\mathbb{C}}

\newcommand{\Sph}[1]{{\mathbb S}^{#1}}  

\newcommand{\del}{\partial}           
\newcommand{\rec}[1]{\tfrac{1}{#1}}   
\newcommand{\iso}{\approx}            
\newcommand{\tens}{\otimes}           

\newcommand{\maps}{\longmapsto}       
\newcommand{\longto}{\longrightarrow} 

\newcommand{\then}{\Longrightarrow}   

\newcommand{\inner}[1]{\langle #1\rangle} 
\newcommand{\Inner}[1]{\bigl\langle #1\bigr\rangle} 
\newcommand{\Rinner}[1]{\inner{#1}\Sub{\aR}}
\newcommand{\RInner}[1]{\Inner{#1}\Sub{\aR}}
\newcommand{\Cinner}[1]{\inner{#1}\Sub{\Ce}}

\newcommand{\norm}[1]{\left\|#1\right\|}  

\newcommand{\rest}[1]{|_{#1}}                          
\newcommand{\srest}[1]{\raisebox{-2pt}{$|_{#1}$}}      

\newcommand{\End}{\operatorname{End}} 

\newcommand{\re}{\operatorname{Re}}   

\newcommand{\supp}{\operatorname{supp}} 
\newcommand{\proj}[1]{\operatorname{proj}\!\mraise{-3pt}{_{#1}} }  

\newcommand{\T}[1]{T_{#1}}             

\newcommand{\TM}{\T{M}}

\newcommand{\uaR}{\underline{\aR}}     
\newcommand{\uCe}{\underline{\Ce}}     

\newcommand{\cont}{\,\lrcorner\,}      
\newcommand{\Nabla}[1]{\nabla\adjustnabla{#1}}

\newcommand{\PDeq}{\stackrel{_{PD}}{=}} 

\newcommand{\Ks}{K^{*}}                

\newcommand{\Ga}{{\mathcal G}}          
\newcommand{\Con}{{\mathcal C}onn}       

\newcommand{\cli}{\,\mraise{3pt}{_{\bullet}}\, }     
\newcommand{\varcli}{\mbox{}_{^{\bullet}}} 

\newcommand{\Spinc}{\ensuremath{Spin^{\Ce}}}   

\newcommand{\pinc}{pin$^{\!{\bf C}}$} 
\newcommand{\spinc}{s\pinc}         

\newcommand{\W}{{\mathcal W}}           

\newcommand{\D}{{\mathcal D}}

\newcommand{\twi}{\goth{s}}         
\newcommand{\Twi}{\twi\Lambda^{+}}  

\newcommand{\GW}{\Gamma(\W^{+})}    
\newcommand{\Wp}{\W^{+}}            
\newcommand{\Lp}{\Lambda^{+}}       

\newcommand{\I}{\goth{z}}

\begin{document}

{\flushleft
appeared in:\\
\smallskip
{\sl Communications in Contemporary Mathematics},\\
vol.~4 (2002), no.~1, pp.~45--64\\
{\small\tt http://ejournals.wspc.com.sg/ccm/ccm.html}
}

\vskip40pt

\title[Nowhere zero harmonic spinors and self-dual $2$-forms]{Nowhere zero harmonic spinors\\ and self-dual $2$-forms}
\author{Alexandru Scorpan}
\address{Department of Mathematics, University of California, Berkeley\\ 970 Evans Hall, Berkeley, CA 94720}
\email{scorpan@math.berkeley.edu}
\urladdr{www.math.berkeley.edu/\textasciitilde scorpan}
\date{June 6, 2001}

\subjclass[2000]{Primary 53C27; Secondary 32Q60, 53D05, 53C55, 57N13}
\keywords{spinor, four-manifold, Dirac operator, almost-complex, symplectic, \kah}

\begin{abstract}

Let $M$ be a closed oriented $4$-manifold, with \riem\ metric $g$, 
and a \spinc-structure induced by an \ac\ structure $\omega$.
Each connection $A$ on the determinant line bundle
induces a unique connection $\nabla^{A}$,
and Dirac operator $\D^{A}$ on spinor fields.
Let $\sigma:\Wp\to\Lp$ be the natural squaring map, taking \SD\ 
(= positive) spinors to \SD\ $2$-forms.

In this paper, 
we characterize the \SD\ $2$-forms that are images of \SD\ spinor fields 
through $\sigma$. They are those $\alpha$ for which (off zeros)
$c_{1}(\alpha)=c_{1}(\omega)$, where $c_{1}(\alpha)$ is a 
suitably defined Chern class.
We also obtain the formula:
$\norm{\phi}^{2} \D^{A}\phi
  =i\bigl( 2\,d^{*}\sigma(\phi) + 
  \Rinner{\nabla^{A}\phi,\,i\phi} \bigr)
  \varcli\phi$.

Using these, we establish a bijective correspondence between:
$\bigl\{$\kah\ forms $\alpha$ 
  compatible with a metric scalar-multiple of $g$, 
  and with $c_{1}(\alpha)=c_{1}(\omega) \bigr\}$
and 
$\bigl\{$gauge classes of pairs $(\phi,A)$,  with 
  $\nabla^{A}\phi= 0 \bigr\}$,
as well as a bijective correspondence between:
$\bigl\{$Symplectic forms $\alpha$ 
  compatible with a metric conformal to $g$, 
  and with $c_{1}(\alpha)=c_{1}(\omega) \bigr\}$
and 
$\bigl\{$gauge classes of pairs $(\phi,A)$, with 
  $\D^{A}\phi= 0$, and 
  $\Rinner{\nabla^{A}\phi,\,i\phi}= 0$, and 
  $\phi$  nowhere-zero$\bigr\}$.
\end{abstract}
\maketitle

Through \cite{witten}, the celebrated \SW\ monopole equations were 
introduced into $4$-dimensional topology, thus 
bringing \spinc\ geometry to the frontline of mathematical research.

For a closed oriented $4$-manifold 
endowed with a fixed \riem\ metric $g$, 
the ingredients of the \SW\ equations are: a \spinc-\str\ 
with spinor bundles $\W^{\pm}$, 
a connection $A$ on their determinant line bundle $L$, an 
associated Dirac operator $\D^{A}:\GW\to\Gamma(\W^{-})$, 
and a quadratic map $\sigma:\Wp\to\Lp$ taking \SD\ 
(= positive%
\footnote{%
Instead of ``\SD\ spinor'', a more customary terminology would be ``positive spinor''. We prefer 
the former, which is used in the classical 
paper \cite{AHS} and seems better suited to the peculiarities of 
dimension $4$.
}%
) 
spinors to \SD\ $2$-forms. 
Solutions to the \SW\ equations are gauge classes of pairs $(\phi,A)$,
made from a spinor field $\phi\in\GW$ and a connection $A$ on $L$.
One of the \SW\ equations is $\D^{A}\phi= 0$, while the second 
equates $\sigma(\phi)$ and part of the curvature of $A$. 

While this paper is not concerned with 
\SW\ theory, its focus is on the same creatures  
mentioned above. 
The paper grew out of an attempt to project the spinorial 
world down to the $4$-manifold. The device mediating this descent was 
to be the squaring map $\sigma$. 

In this paper we study the images 
$\sigma(\phi)$ of gauge classes of pairs $(\phi,A)$ satisfying the 
equation $\D^{A}\phi= 0$. 
First, we determine which \SD\ $2$-forms are images of \SD\ spinor 
fields through $\sigma$. 
Then we determine when these $2$-forms are 
\kah\ or symplectic. 
The main results are summarized in Table~1,
and proved as Theorems \ref{thm-a} and \ref{thm-b}.
(Note that we restrict ourselves to the special case of
\spinc-\str s that come from \emph{\ac} \str s.)
An excellent general reference for related 
material (genuine spin-\str s) is \cite{spingeom}. For the \spinc\ point 
of view (and \SW\ theory), see \cite{morgan}.

%
\begin{table}[htb]
\caption{There is a bijection, 
  established through $\alpha=\sigma(\phi)$, 
  between the left and the right sides
  (when $H^{2}(M;\Zi)$ has no $2$-torsion).}

\newcommand{\aspace}{\mbox{}\quad }

\begin{center}
\small
\begin{tabular}{|p{177pt}|p{158pt}|}
\multicolumn{2}{c}{}\cr
\hline

&\cr
\ \kah\ forms $\alpha$ &
\ Gauge classes of pairs $(\phi,A)$ \cr
\aspace compatible with a metric multiple of $g$  &
\aspace with $\nabla^{A}\phi= 0$ \cr
\aspace and with $c_{1}(\alpha)=c_{1}(M)$.&
\aspace and $\phi$ not constantly-zero.\cr
&\cr
\hline
&\cr

\ Symplectic forms $\alpha$ &
\ Gauge classes of pairs $(\phi,A)$ \cr
\aspace compatible with a metric conformal to $g$ &
\aspace with $\D^{A}\phi= 0$ and $\RInner{\nabla^{A}\phi,\ i\phi}= 0$ \cr
\aspace and with $c_{1}(\alpha)=c_{1}(M)$. &
\aspace and $\phi$ nowhere-zero.\cr
&\cr
\hline
&\cr

\ Self-dual nowhere-zero $2$-forms $\alpha$ &
\ Gauge classes of pairs $(\phi,A)$ \cr
\aspace with $c_{1}(\alpha)=c_{1}(M)$ &
\aspace with $\D^{A}\phi= 0$, and $\phi$ nowhere-zero. \cr
&\cr

\hline
\end{tabular}
\end{center}

\label{table}
\end{table}
%

Section~1 presents \ac\ \str s  as constant-length sections of $\Lp$, 
then relates Chern classes of \ac\ \str s with intersection 
theory in $\Lp$. A main result is:
\begin{equation}
  c_{1}(\alpha)+c_{1}(\beta) \PDeq 2\cdot\alpha\cap\beta\rest{M} 
\tag{Thm.~\ref{thm-chern}}
\end{equation}
where on the right side we took intersections in the sphere bundle 
of $\Lp$ and projected them as  oriented surfaces in $M$.

Section~2 first reviews \spinc-\str s associated to an \ac\ 
\str\ $\omega$,
and defines the squaring map $\sigma:\Wp\to\Lp$. 
Then it attacks the lifting problem for $\sigma$, that is, the 
problem of which \SD\ $2$-forms $\alpha$ are 
images $\sigma(\phi)$ of spinor fields. The main result is that 
$\alpha$ can be lifted \Iff, off the zeros of $\alpha$, 
the class $c_{1}(\alpha)$ coincides with the Chern class 
$c_{1}(\omega)$ of the \ac\ \str\ $\omega$ 
that generated the \spinc-\str\ (Thms.~\ref{thm-lift} 
and \ref{cor-lift}).

Section~3 discusses connections on $\W^{\pm}$ and defines the 
Dirac operator $\D^{A}$. 
Then it works at separating the contributions of the \LC\ connection and of $A$
in building the Dirac operator. We obtain the formula
\begin{equation}
\norm{\phi}^{2} \D^{A}\phi
  =i\bigl( 2\,d^{*}\sigma(\phi) + \Rinner{\nabla^{A}\phi,\,i\phi} \bigr)
  \cli\phi   
\tag{Thm.~\ref{thm-split.dirac}}
\end{equation}

Section~4 gathers together all the results obtained along the way, 
and proves the statements from Table~\ref{table}. 

\bigskip %

While this work was being completed, 
the paper \cite{4koreans} was electronically published. 
It contains two partial results 
in the direction of our work
(see Remarks \ref{rk-4k.1} and \ref{rk-4k.2} for a discussion).
Developed independently, 
our paper is strengthening their results.


\section{Almost-Complex Structures}
\label{sec-ac}

Let $M$ be a closed oriented $4$-manifold, endowed with a \emph{fixed} \riem\ 
metric $g$, and its \LC\ connection $\nabla$. The same notation $\nabla$ 
will denote the connections induced on the tensor bundles of $M$.
Also, using the metric, we will systematically identify $\TM$ and $\TM^{*}$. 
Throughout the paper, an ``$x$'' will denote a generic point of $M$.

Let $\Lp$ denote the $3$-plane bundle of \SD\ $2$-forms on $M$, 
and $\Gamma(\Lp)$ denote the space of \SD\ $2$-forms (sections in $\Lp$).
Also, let $\Twi$ denote the bundle of $2$-spheres of radius $\sqrt{2}$ in $\Lp$. 
(The latter is also known as the \defemph{twistor bundle} of $M$.)
At a point $x$, each $\alpha\in\Twi\rest{x}$ can be 
written as $\alpha=e_{1}\wedge e_{2}+e_{3}\wedge e_{4}$ for some 
orthonormal orienting basis $\{e_{k}\}$ in $\TM\rest{x}$. Therefore it 
defines an \ac\ structure $J_{\alpha}$ 
on $\TM\rest{x}$ by $J_{\alpha}(e_{1})=e_{2}$ 
and $J_{\alpha}(e_{3})=e_{4}$. The \ac\ structure $J_{\alpha}$ is 
orthogonal for $g\srest{x}$, \ie $g(J_{\alpha}v,\,J_{\alpha}w)=g(v,w)$.
Conversely, given an orthogonal \ac\ structure $J$ 
on $\TM\rest{x}$, there is an associated $2$-form 
$\alpha\Sub{J}\in\Twi\rest{x}$ given by 
$\alpha\Sub{J}\srest{x}=e_{1}\wedge e_{2}+e_{3}\wedge e_{4}$
for any orthonormal orienting basis $\{e_{k}\}$ with $J(e_{1})=e_{2}$ and 
$J(e_{3})=e_{4}$. (In general, a $2$-form $\alpha\in\Twi\rest{x}$ and 
an orthogonal \ac\ structure $J$ on $\TM\rest{x}$ are associated through
$\alpha(v,w)=g(Jv,w)$.
From now on, \emph{all} \ac\ structures 
considered will be orthogonal.) 

In conclusion, the bundle $\Twi$ can be seen as the bundle 
of orthogonal \ac\ structures on fibers of $\TM$. 
Global \ac\ structures on $M$ can be identified with global sections of 
$\Twi$. 

An \ac\ structure $\alpha$ on $M$ makes $\TM$ a complex 
bundle. The complex-line bundle $\Ks_{\alpha}=\det_{\Ce}(\TM,\alpha)$ 
is called the \defemph{anti-canonical bundle} of $\alpha$. Its Chern 
class $c_{1}(\Ks_{\alpha})=c_{1}(\TM,\alpha)\in H^{2}(M;\Zi)$ is called the 
\defemph{Chern class} of the \ac\ structure $\alpha$ and will be 
denoted by $c_{1}(\alpha)$. 
The bundle $\Ks_{\alpha}$ can also be 
recovered as the oriented complement of $\aR\alpha$ in 
$\Lambda^{+}$, because of the oriented splitting 
\[ \Lp=\aR\alpha\oplus\Ks_{\alpha} \]

Each nowhere-zero \SD\ $2$-form $\beta\in\Gamma(\Lp)$ has a unique 
projection $\twi(\beta)$ in $\Gamma(\Twi)$ and thus defines an \ac\ 
structure with a well-defined Chern class $c_{1}(\beta)$.

\bigskip %

In what follows we will study  the relations between \ac\ structures, 
as sections in $\Twi$, and their Chern classes.

Let $\alpha, \beta\in\Gamma(\Twi)$ be two sections, seen as 
oriented $4$-submanifolds in $\Twi$. Make one transverse 
to the other. Their intersection $\alpha\cap\beta$ 
will be an \emph{oriented} surface in $\Twi$ that projects nicely to 
an oriented surface in $M$.
The latter will be denoted by by $\alpha\cap\beta\rest{M}$. 

Looking at the signs of the intersections, one can see that
we have: $\alpha\cap\beta=\beta\cap\alpha$, and also: 
$\alpha\cap\beta\rest{M}=- \bigl( (-\alpha)\cap(-\beta)\rest{M} 
\bigr)$, where the minus represents reversal of orientation.

\begin{figure}[htb]
\begin{center}
\begin{picture}(240,120)(-20,0)

\put(-20,26.6){\line(1,0){240}}
\put(-20,26.6){\line(1,2){40}}

\put(100,60){\bigcircle{120}}


\qbezier(40,60)(40,40)(100,40)
\qbezier(100,40)(160,40)(160,60)
\qbezier[20](40,60)(40,80)(100,80)
\qbezier[20](100,80)(160,80)(160,60)

\put(100,120){\makebox(0,0){$\bullet$}}
\put(100,0){\makebox(0,0){$\bullet$}}

\put(96,109){\large$\alpha$}
\put(87,5.5){\large$-\alpha$}

\put(140,110){\large$\Twi\srest{x}$}
\put(0,38){\large$\Ks_{\alpha}\srest{x}$}
\put(-20,100){\large$\Lambda^{+}\srest{x}$}

\put(100,60){\makebox(0,0){\scriptsize$\bullet$}}

\put(200,120){\vector(0,-1){60}}
\put(203.3,86.6){{\footnotesize $proj$}}
\put(200,0){\vector(0,1){20}}

\end{picture}
\end{center}
\caption{Almost-complex structure and bundles}
\label{fig-proj}
\end{figure}
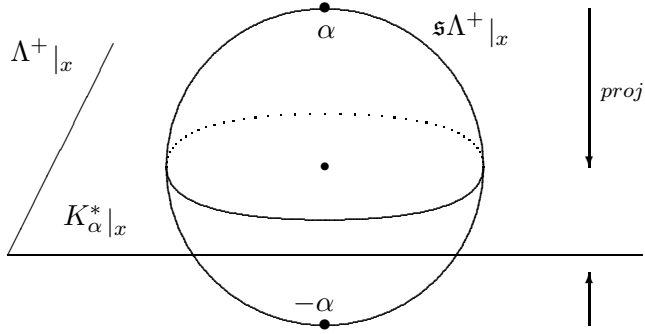

Let $\alpha, \beta\in\Gamma(\Twi)$ be two \ac\ structures, with 
anti-canonical bundles $\Ks_{\alpha}$ and $\Ks_{\beta}$. 
Make $\beta$ transverse to $\alpha$ in $\Twi$, then project it on 
$\Ks_{\alpha}$ through the composition 
$\Twi\subset\Lambda^{+}=\aR\alpha\oplus\Ks_{\alpha}
  \stackrel{pr}{\longto}\Ks_{\alpha}$
(see \Figref{fig-proj}). 
Call $\beta_{K}$ the resulting section of $\Ks_{\alpha}$.
%
%
Since $\beta$ is transverse to $\alpha$ in $\Twi$, 
it follows that the projection $\beta_{K}$
is transverse to the zero section of $\Ks_{\alpha}$. 
The  intersection points of $\beta_{K}$ with the zero section in 
$\Ks_{\alpha}$
correspond to the  intersection points of $\beta$ with $\alpha$ 
or with $-\alpha$ in $\Twi$. 
More precisely, after matching signs:
\[ \beta_{K}\cap 0\,\rest{M}
  =\beta\cap\alpha\rest{M}-\beta\cap(-\alpha)\rest{M} \]
But the zero locus of a generic section in $\Ks_{\alpha}$ is the \PD\ 
of the Chern class $c_{1}(\alpha)$. We obtained:

\begin{lemma}%
For all sections $\alpha,\beta$ in $\Twi$, we have: 
\[ c_{1}(\alpha)\PDeq 
  \beta\cap\alpha\rest{M}-\beta\cap(-\alpha)\rest{M} \]
In particular: 
$c_{1}(\alpha)\PDeq  \alpha\cap\alpha\srest{M}$.
\end{lemma}
But 
$\beta\cap(-\alpha)\rest{M}=-\bigl( (-\beta)\cap\alpha\rest{M} \bigr)$, 
so we can write as well:
\[ c_{1}(\alpha)\PDeq 
  \alpha\cap\beta\rest{M}+\alpha\cap(-\beta)\rest{M} \]
On the other hand, reversing the r\^oles of $\alpha$ and $\beta$ in 
the Lemma:
\[ c_{1}(\beta)\PDeq 
  \alpha\cap\beta\rest{M}-\alpha\cap(-\beta)\rest{M} \]
Adding, we get:

\begin{theorem}%
\label{thm-chern}
For all sections $\alpha, \beta\in\Gamma(\Twi)$, we have
\[ c_{1}(\alpha)+c_{1}(\beta) \PDeq 2\cdot \alpha\cap\beta\rest{M} \]
\end{theorem}
A consequence is that a cocycle for the cohomology class 
$c_{1}(\alpha)+c_{1}(\beta)$ 
can be written as:
$s^{2}\maps 2(\alpha\cap\beta\rest{M}) \cap s^{2}$, 
for all singular $2$-simplices $s^{2}$ in $M$ transverse to the surface 
$\alpha\cap\beta\rest{M}$.

Without trouble, we can restate this as:

\begin{corollary}
\label{cor-cocycle}
A cocycle for $c_{1}(\alpha)+c_{1}(\beta)$
can be written as:
\[ s^{2} \maps 2\,\alpha(s^{2})\cap\beta 
  \qquad\text{in }\Twi \]
for all singular $2$-simplices $s^{2}$ in $M$ 
\st\ their image $\alpha(s^{2})$ is
transverse to the $4$-submanifold $\beta$ in $\Twi$ .
\end{corollary}%

\bigskip%

We conclude this section with a few definitions:

An \ac\ structure $\alpha\in\Gamma(\Twi)$ 
is called a \defemph{\kah\ form} compatible with the metric $g$
\Iff\ $\nabla\alpha= 0$,
where $\nabla$ is the connection on $\Lp$ induced from $g$'s \LC\ connection. 
Indeed, in that case $J_{\alpha}$ is integrable, and 
$(M,g,J_{\alpha})$ is \kah, with \kah\ form $\alpha$.

Let $\alpha$ be a nowhere-zero \SD\ $2$-form with
$\nabla\alpha= 0$. Then $\alpha$ has constant length. 
Modifying the metric from $g$ to $g'$ 
by multiplying with a suitable constant, we can make
$\alpha$ of length $\sqrt{2}$.
Thus $\alpha$ becomes a section in $(\Twi)'$ and 
determines an \ac\ structure $J_{\alpha}'$.
Since the \LC\ connection stays the same, $\nabla'=\nabla$, 
we conclude that $(M,g',J_{\alpha}')$ is \kah.
In conclusion:

\begin{lemma}%
\label{lemma-kah}
A nowhere-zero \SD\ $2$-form $\alpha\in\Gamma(\Lp)$ has $\nabla\alpha= 0$ 
\Iff\ $\alpha$  is a \kah\ form compatible with a metric scalar-multiple of $g$.
\end{lemma}%

A non-degenerate $2$-form $\alpha$ with $d\alpha=0$ is called 
\defemph{symplectic}. 
If a symplectic form $\alpha$ is \SD\ and of constant length $\sqrt{2}$, 
then $\alpha$ is also called \emph{compatible} with the metric $g$.
(And $(M,g,\alpha)$ is called an \defemph{almost-\kah} manifold.)

\begin{lemma}%
\label{lemma-symp}
A nowhere-zero \SD\ $2$-form $\alpha\in\Gamma(\Lp)$ has $d\alpha=0$
\Iff\ $\alpha$ is a symplectic form compatible with a metric 
conformal to $g$.
\end{lemma}
\begin{proof}
A first thing to notice is that 
``nowhere-zero and \SD'' insures ``non-degenerate''.
When suitably rescaling the metric, one always preserves the bundle $\Lp$, 
but can make $\norm{\alpha}\equiv\sqrt{2}$.
\end{proof}%
 
One goal of this paper is to investigate how certain  spinor 
fields determine such \kah\ or symplectic structures on $M$.

%

\section{The Lifting Problem}
\label{sec-lift}

\subsection{S\pinc-structures}

Since our manifold is oriented \riem, the bundle $\TM$ has structure 
group $SO(4)= SU(2)^{+}\times SU(2)^{-}\,\big/\pm1$.
The complex spin group can be described as
$\Spinc(4)= \Sph{1}\times SU(2)^{+}\times SU(2)^{-}\,\big/\pm1$,
and projects canonically to $SO(4)$.
The manifold $M$ is said to admit a \defemph{\spinc-structure} \Iff\ the 
$SO(4)$-cocycle of $\TM$ can be lifted to a $\Spinc(4)$-cocycle.
On $4$-manifolds, \spinc-structures always exist.

Suppose given an \ac\ structure $\omega\in\Gamma(\Twi)$. Then the 
$SO(4)$-cocycle of $\TM$ can be reduced to an $U(2)$-cocycle.
More, since 
$U(2)= \Sph{1}\times SU(2)^{-}\,\big/\pm1
  \ \subset\  SU(2)^{+}\times SU(2)^{-}\,\big/\pm1=SO(4)$,
and 
$\Spinc(4)= \Sph{1}\times SU(2)^{+}\times SU(2)^{-}\,\big/\pm1$,
there is an obvious canonical map $U(2)\to\Spinc(4)$ 
covering the inclusion $U(2)\subset SO(4)$.
Then one can use this map to lift the $U(2)$-cocycle of $\TM$ to a 
$\Spinc(4)$-cocycle. 
Therefore, each \ac\ structure induces a canonical 
\spinc-structure on $M$. 

Choose some \spinc-structure on $M$, \ie a $\Spinc(4)$-cocycle for $\TM$.
Then the  natural projections
$\Spinc(4)\to U(2)^{\pm}= \Sph{1}\times SU(2)^{\pm}\,\big/\pm1$
induce $U(2)$-cocycles for hermitian complex-plane bundles $\W^{\pm}$, 
called the \defemph{bundles of \SD/\ASD\ spinors}. 
The complex-line bundle $L=\det_{\Ce}\W^{\pm}$ is called the
\defemph{determinant bundle} of the \spinc-structure.
If the \spinc-\str\ comes from an \ac\ \str\ $\omega$, then 
its determinant bundle $L$ is exactly the 
anti-canonical bundle $\Ks_{\omega}$ of the \ac\ structure.

A \spinc-structure also determines a \defemph{Clifford 
multiplication} map $\TM\times\W^{+}\to\W^{-}$, denoted by 
$(v,\phi)\maps v\cli\phi$. It is modelled on quaternionic 
multiplication.
In particular, if $v\cli\phi=0$, then $v=0$ or $\phi=0$.
Clifford multiplication induces an action $\Lp\times\Wp\to\Wp$.
When the \spinc-\str\ is associated to an \ac\ \str\ $\omega$,
the Clifford action of the $2$-form $\omega$ splits $\Wp$ into its 
$\mp 2i$-eigenbundles as 
\[ \Wp=\uCe\oplus\Ks_{\omega} \]
On the \ASD\ side, we have that $\W^{-}\iso (\TM,\omega)$ 
(as complex bundles).

Let $\End_{0}(\Wp)$ denote the space of 
\emph{traceless} $\Ce$-endomorphisms of $\Wp$.
The Clifford action $\Lp\times\Wp\to\Wp$ identifies
$\End_{0}(\Wp)\iso\Lp\tens\Ce$. 
For every $\phi\in\Wp\rest{x}$, 
we have that $\phi\tens\phi^{*}$ is an endomorphism of 
$\Wp$ with traceless part 
$\phi\tens\phi^{*}-\rec{2}\norm{\phi}^{2}id$. The latter corresponds to a 
purely imaginary $2$-form $i\sigma(\phi)\in i\Lp\rest{x}$.
This defines the \defemph{squaring map} 
\[ \sigma:\Wp\to\Lp \]
It is uniquely characterized by: 
\[ \sigma(\phi)\cli\phi=-i\tfrac{\norm{\phi}^{2}}{2}\phi \]
If the \spinc-\str\ comes from an \ac\ \str\ $\omega$, then: 
\[\begin{CD}
\phi &\,\in\,&\uCe &\subset\Wp&  \ \mbox{of length 2}\ & 
  \quad\then\qquad 
  &\omega\cli\phi=-2i\cdot\phi&
  \quad\then\qquad
  &\sigma(\phi)=+\omega \\
\phi &\,\in\,&\Ks &\subset\Wp&   \ \mbox{of length 2}\ &
  \quad\then\qquad 
  &\omega\cli\phi=+2i\cdot\phi&
  \quad\then\qquad
  &\sigma(\phi)=-\omega 
\end{CD}\]

The squaring map has the property that $\sigma(e^{if}\phi)=\sigma(\phi)$. 
Thus, if we restrict $\sigma$ to a sphere in a fiber $\Wp\rest{x}$, 
then it behaves like a Hopf map $\Sph{3}\to \Sph{2}$.
And since $\norm{\sigma(\phi)}=\rec{2\sqrt{2}}\norm{\phi}^{2}$, one 
can further look at $\sigma$ fiber-wise
as a sort of squared-cone on the Hopf map. More precisely, if 
${\goth h}: \Sph{3}\to \Sph{2}$
is the Hopf map between 
the sphere (of radius 2) in $\Wp\rest{x}$ 
and the sphere (of radius $\sqrt{2}$) in $\Lp\rest{x}$, 
then 
\[ \sigma(a\phi)=a^{2}\cdot {\goth h}(\phi) \]
for all $a\in\aR$ and $\phi\in \Sph{3}$.

\bigskip %

\begin{remark}%
\label{rk-id}
A generic \spinc-structure has the following property: $\Wp$ admits a 
nowhere-zero section \Iff\ the \spinc-structure comes from some 
\ac\ structure on $M$. 
More, in this case there is a \defemph{distinguished spinor field} 
$\I\in\GW$ of constant-length, which is uniquely
characterized by the property that, for any $v\in \TM$, 
we have $v\cli\I=v$ (equality seen through $\W^{-}\iso \TM$). 
In particular, in the splitting $\Wp=\uCe\oplus\Ks$ we have 
$\uCe=\Ce\cdot\I$, and also $\sigma(\I)=\rec{4}\,\omega$.
\end{remark}%

From now on, we \emph{fix} an \ac\ structure $\omega\in\Gamma(\Twi)$ and its 
associated \spinc-structure, with spinor bundles $\W^{\pm}$, and 
determinant bundle $\Ks=\Ks_{\omega}$.

The squaring map is surjective fiber-wise, but not surjective on 
sections: There are section in $\Lp$ that are not images through 
$\sigma$ of sections of $\Wp$. 
This \emph{lifting problem} for $\sigma$ will 
be investigated in the rest of this section. 
The topological obstruction will turn out to be $c_{1}(\alpha)-c_{1}(\omega)$.

\subsection{The setting}

\newcommand{\TW}{\twi\W^{+}} 

For every vector bundle $E\to B$, we will denote by 
$\twi E\to B$ a sphere bundle (of some chosen radius) of $E$.
Since the lifting problem is topological in nature, 
we will blissfully ignore for a while all 
positive scalars involved in the expressions for the squaring map 
$\sigma$ (such as $\tfrac{\norm{\phi}^{2}}{2\sqrt{2}}$, $\sqrt{2}$, etc.), 
and consider 
$\sigma$ simply as a fiber-wise Hopf map from the $3$-sphere bundle 
$\TW$ to the $2$-sphere bundle $\Twi$.

Let $\alpha:M\to\Twi$ be a section. We try to 
lift $\alpha$ to a section $\phi:M\to\TW$ \st\ $\sigma(\phi)=\alpha$.

First note that $\sigma:\TW\to\Twi$ is a circle-bundle over 
$\Twi$, made from patching Hopf fibrations $ \Sph{3}\to \Sph{2}$.
Therefore it has a well-defined 
Chern class $c_{1}(\sigma)=c_{1}(\TW\to\Twi)\ \in H^{2}(\Twi;\Zi)$.

Pull-back the bundle $\sigma:\TW\to\Twi$ through $\alpha:M\to\Twi$ 
to get the circle-bundle $\alpha^{*}(\TW)\to M$, fitting in the diagram:
\[\begin{CD}
\alpha^{*}(\TW) @>>> \TW\\
@VVV @VV\sigma V\\
M @>>\alpha >\Twi
\end{CD}\]
A lift of $\alpha:M\to\Twi$ to a section $\phi:M\to\TW$ is equivalent 
to a section in the circle-bundle $\alpha^{*}(\TW)\to M$. 
But the latter bundle has a section \Iff\ it is trivial. 
That means: \Iff\ its Chern class is zero. 
But 
\[ c_{1}\bigl( \alpha^{*}(\TW)\to M \bigr) 
  =\alpha^{*}[ c_{1}(\sigma)] \]
where the $\alpha^{*}$ on the right is the induced morphism 
$\alpha^{*}:H^{2}(\Twi;\Zi)\to H^{2}(M;\Zi)$.
Therefore:

\begin{lemma}%
\label{lemma-lift}
A nowhere-zero \SD\ $2$-form 
$\alpha$  can  be lifted to a spinor field \Iff\ 
$\alpha^{*}[ c_{1}(\sigma)]=0$ in $H^{2}(M;\Zi)$.
\end{lemma}

\begin{remark}%
For example, take the two distinguished $2$-forms 
$\omega$ and $-\omega$. 
Since $\W^{+}=\uCe\oplus\Ks$ and $\sigma(\twi\uCe)=\omega$ and 
$\sigma(\twi\Ks)=-\omega$, 
we must have the pull-back diagrams:
\[\begin{CD}
  \twi\uCe @>>> \TW\\ @VVV @VV\sigma V\\ M @>>\omega> \Twi
  \end{CD}
\qquad\qquad\qquad 
\begin{CD}
  \twi\Ks @>>> \TW\\ @VVV @VV\sigma V \\ M @>>-\omega > \Twi
\end{CD}\]
and therefore:
$\omega^{*}\bigl( c_{1}(\sigma) \bigr)=0$
and
$(-\omega)^{*}\bigl( c_{1}(\sigma) \bigr)=c_{1}(\Ks)$.
And indeed, $\omega$ lifts immediately to the distinguished spinor 
$\I$ (see \ref{rk-id}), while in general
$-\omega$ cannot be lifted at all. 

In what follows, we will first identify $c_{1}(\sigma)$ as the \PD\ 
of $-\omega$ in $\Twi$, and then show that 
$2\,\alpha^{*}[c_{1}(\sigma)]=c_{1}(\alpha)-c_{1}(\omega)$.
\end{remark}%

\subsection{A section in $\sigma$}

\begin{figure}[htb]
\begin{center}
\begin{picture}(360,180)(0,0)

\put(315,90){\bigcircle{90}}


\qbezier[25](270,90)(270,80)(315,80)
\qbezier[25](315,80)(360,80)(360,90)
\qbezier[10](270,90)(270,100)(315,100)
\qbezier[10](315,100)(360,100)(360,90)

\put(315,135){\makebox(0,0){\scriptsize$\bullet$}}
\put(315,45){\makebox(0,0){\scriptsize$\bullet$}}
\put(312,140){$\omega$}
\put(304,35){$-\omega$}
\put(330,170){\large$\Twi\srest{x}$}

\put(195,90){\vector(1,0){60}}
\put(220,97){\large$\sigma$}

\qbezier(0,180)(90,180)(90,90)
\qbezier(90,90)(90,70)(85.5,55)
\qbezier(0,0)(65,0)(83,47)

\put(70,160){$\sigma^{-1}(\omega)$}
\put(90,90){\makebox(0,0){\scriptsize$\bullet$}}
\put(95,87){\large$\I$}

\qbezier(90,130)(180,130)(180,90)
\qbezier(180,90)(180,50)(90,50)
\qbezier(90,50)(0,50)(0,90)
\qbezier(0,90)(0,130)(80,130)
\put(130,42){$\sigma^{-1}(-\omega)$}

\put(130,170){\large$\TW\srest{x}$}

\end{picture}
\end{center}
\caption{A piece of $\sigma$ as a Hopf fibration}
\label{fig-hopf.sigma}
\end{figure}
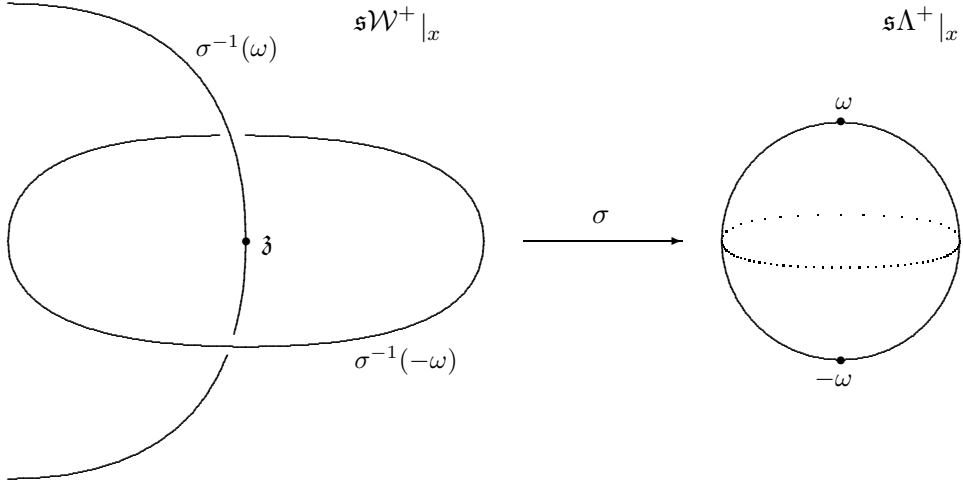

Think of $\sigma:\TW\to\Twi$ as a bundle made up from patching Hopf
fibrations over each $2$-sphere fiber in $\Twi$
(see \Figref{fig-hopf.sigma}).
Since $\Lambda^{+}=\aR\omega\oplus\Ks$,  we can 
imagine that, as we move among
the  $2$-spheres of $\Twi$, the poles $\omega$ and $-\omega$
stay still, while the equator moves around and twists.
Since $\W^{+}=\uCe\oplus\Ks$ (with $\uCe=\Ce\cdot\I$),
we can also imagine that, 
as we move among the $3$-spheres of $\TW$, the point
determined by the spinor field $\I$ (see \ref{rk-id})
stays still as well.
The circle $ \Sph{1}\cdot\I_{x}=\twi\uCe\,\rest{x}$ of $\TW\rest{x}$
is sent through $\sigma$ to
the $\omega_{x}$-pole in $\Twi\rest{x}$.
The circle $\twi\Ks\rest{x}$ is sent through $\sigma$ to $-\omega_{x}$.

To determine the Chern class of the circle-bundle 
$\sigma:\TW\to\Twi$, we will build a section $s$ of its associated 
complex-line bundle ${\cal L}\to\Twi$,
by patching sections over each $2$-sphere of $\Twi$. 

Lift the pole $\omega_{x}$ of the $2$-sphere $\Twi\rest{x}$
to $\I_{x}$, and thus define $s(\omega)=\I$.
Further, 
since $\I_{x}$ and the circle $\sigma^{-1}(-\omega_{x})$ 
over $-\omega_{x}$ are well-determined in each 
$3$-sphere, we can choose \emph{uniformly}
a disk $D_{x}$ in $\TW\rest{x}$,
centered at $\I_{x}$ and bounded by the circle 
$\sigma^{-1}(-\omega_{x})$. 
  
Over each point $\lambda$ in the $2$-sphere $\Twi\rest{x}$
lies a circle $\sigma^{-1}\lambda$ in the $3$-sphere $\TW\rest{x}$. 
If we chose the disk $D_{x}$ nicely enough,
the circle $\sigma^{-1}\lambda$ cuts $D_{x}$ in exactly
\emph{one} point. 
Define $s(\lambda)$ to be that intersection point: 
$s(\lambda)=\sigma^{-1}\lambda\cap D_{x}$.

The only exception is when $\lambda=-\omega_{x}$. 
The circle over $-\omega_{x}$ lies entirely in $D_{x}$ as its boundary. 
Define $s(-\omega_{x})=0$. To make the resulting $s$ 
continuous, modify the lengths of all 
other  $s(\lambda)$'s, so that, as $\lambda$ approaches $-\omega_{x}$ 
in $\Twi\rest{x}$, $s(\lambda)$ will approach $0$.

This defines $s$ in one Hopf patch of $\sigma$. 
We repeat the same procedure over the whole $\Twi$.
With some care, in the end
we have a well-defined continuous 
section $s:\Twi\to {\cal L}$ that has 
zeros exactly along the image of $-\omega:M\to\Twi$. That means:

\begin{lemma}%
The Chern class $c_{1}(\sigma)\in H^{2}(\Twi;\Zi)$ 
of the bundle $\sigma:\TW\to\Twi$ is the \PD\ of 
$(-\omega)_{*}[M]\in H_{4}(\Twi;\Zi)$. Or, in short:
\[ c_{1}(\sigma)\PDeq  -\omega  \]
\end{lemma}

\subsection{\mbox{}\ldots and lifting}

Since $ c_{1}(\sigma)\PDeq  -\omega $, that means that  
$c_{1}(\sigma)\in H^{2}(\Twi;\Zi)$ can be represented by the cocycle:
$s^{2} \maps s^{2}\cap(-\omega)$,
for all singular $2$-simplices $s^{2}$ in $\Twi$ 
that are transverse to $-\omega$.

For a  section $\alpha:M\to\Twi$, the class 
$\alpha^{*}[c_{1}(\sigma)]\in H^{2}(M;\Zi)$ 
is then represented by the cocycle:
\[ s^{2} \maps \alpha(s^{2})\cap(-\omega) \qquad\text{in }\Twi \]
for all singular $2$-simplices $s^{2}$ in $M$ whose images through $\alpha$ 
are transverse to the $4$-submanifold $-\omega$ of $\Twi$.

But remember from \Corref{cor-cocycle} that 
this is exactly half a cocycle that represents the class 
\[ c_{1}(\alpha)+c_{1}(-\omega) \]
Noticing that $c_{1}(-\omega)=-c_{1}(\omega)$, we have:

\begin{lemma}%
For every section $\alpha:M\to\Twi$, we have
\[ 2\,\alpha^{*}[c_{1}(\sigma)]=c_{1}(\alpha)-c_{1}(\omega) \]
\end{lemma}%

Remember now that 
a section $\alpha$ in the sphere bundle $\Twi$ 
can be lifted to $\TW$
\Iff\ $\alpha^{*}[c_{1}(\sigma)]=0$
(\Lemmaref{lemma-lift}).
That implies $2\,\alpha^{*}[c_{1}(\sigma)]=0$, and so 
$c_{1}(\alpha)=c_{1}(\omega)$. 
Therefore: 

\begin{theorem}%
\label{thm-lift}
Let $\alpha\in\Gamma(\Twi)$.
If $\alpha$ can be lifted to a spinor field, 
then $c_{1}(\alpha)=c_{1}(\omega)$.
If $H^{2}(M;\Zi)$ has no $2$-torsion, 
then the converse is also true.
\end{theorem}
\begin{warning}%
\label{war-torsion}
{\sf For the simplicity of later statements, 
from now on we will assume
that \defemph{$H^{2}(M;\Zi)$ has no $2$-torsion}.}
{\rm Under this assumption, the statement above would be restated as:}
``$\alpha$ can be lifted \Iff\ $c_{1}(\alpha)=c_{1}(\omega)$''.
\end{warning}%

Passing from section in $\Twi$ to nowhere-zero sections in $\Lambda^{+}$ 
is immediate.
Each nowhere-zero section $\alpha:M\to\Lambda^{+}$ has a well-defined 
projection $\twi(\alpha):M\to\Twi$. Then $\alpha$ can be lifted 
\Iff\ $\twi(\alpha)$ can be lifted, using a formula such as
$\sigma(a\phi)=a^{2}\,\sigma(\phi)$ to adjust the lengths.
And Chern classes are, by definition, 
unchanged by $\alpha\maps\twi(\alpha)$.

\begin{corollary}%
A nowhere-zero \SD\ $2$-form $\alpha$ lifts to a spinor field \Iff\ 
$c_{1}(\alpha)=c_{1}(\omega)$.
\end{corollary}%

If $\alpha:M\to\Lambda^{+}$ has zeros, 
repeat the entire discussion on $M\setminus\{\text{zeros}\}$.
We get that $\alpha\srest{\text{off zeros}}$ can be lifted to some 
$\phi\rest{\text{off zeros}}$ \Iff\ 
$c_{1}\bigl(\alpha\srest{\text{off zeros}}\bigr)
 =c_{1}(\omega)\rest{\text{off zeros}}$.
But $\alpha\srest{\text{off zeros}}$ can be lifted \Iff\ $\alpha$ 
can be lifted: Indeed, since $\sigma(0)=0$ and 
$\norm{\sigma(\phi)}=\rec{2\sqrt{2}}\norm{\phi}^{2}$, 
the partial lift $\phi\rest{\text{off zeros}}$ 
can be completed with zeros to get a lift of the whole $\alpha$.
We proved:

\begin{corollary}%
\label{cor-lift}
A \SD\ $2$-form $\alpha$ can be lifted to a continuous spinor field
\Iff\
$c_{1}\bigl( \alpha\srest{\text{off zeros}} \bigr)
  =c_{1}(\omega)\rest{\text{off zeros}}$
in $H^{2}(M\setminus\{zeros\};\Zi)$ 
(assuming $H^{2}(M\setminus\{zeros\};\Zi)$ has no $2$-torsion). 
\end{corollary}
Note that this procedure creates only a \emph{continuous} spinor 
field, which might be not smooth at zeros.

\section{The Splitting Formula}
\label{sec-split}

\subsection{Connections and Dirac operators}

Pick some \spinc-\str\ with spinor bundles $\W^{\pm}$ and 
with determinant bundle $L$. 
Then the \LC\ connection $\nabla$ on $\TM$ and a choice of unitary 
connection $A$ on the determinant line bundle $L$ will determine a 
unique connection $\nabla^{A}$ on the spinor bundles $\W^{\pm}$.

Let $\Con(L)$ denote the space of unitary connections on $L$. It is an 
affine space with model $i\Gamma(\TM^{*})$: any two connections differ 
by an imaginary $1$-form. If $A=A_{0}+2i\theta$, then the corresponding 
spinorial connections have $\nabla^{A}=\nabla^{A_{0}}+i\theta$.

All spinorial connections $\nabla^{A}$ on $\W^{\pm}$
are compatible with the \herm\ metrics:
$\del_{X}\inner{\phi,\psi}
  =\inner{\nabla^{A}_{X}\phi,\psi}+\inner{\phi,\,\nabla^{A}_{X}\psi}$.
The connections $\nabla^{A}$ 
are also compatible with Clifford multiplication:
\[ \nabla^{A}(v\cli\phi)=(\nabla v)\cli\phi+v\cli(\nabla^{A}\phi) \]
for every $v\in\Gamma(\TM)$ and $\phi\in\Gamma(\W^{\pm})$. A similar 
relation holds for Clifford multiplication with $2$-forms.

The \defemph{Dirac operator} $\D^{A}$ is defined by the composition:
\[ \D^{A}: \Gamma(\Wp) \stackrel{\nabla^{A}}{\longto}
  \Gamma(\Wp\tens \TM^{*}) \stackrel{\cli}{\longto} \Gamma(\W^{-}) \]
or, equivalently, by the formula:
\[ \D^{A}\phi=\sum e_{k}\cli\nabla^{A}_{e_{k}}\phi \]
for any local orienting orthonormal frame $\{e_{k}\}$ in $\TM$.
A spinor field $\phi\in\GW$ \st\ $\D^{A}\phi= 0$ is called 
\defemph{harmonic} \emph{for $A$}.

Since, for functions $f:M\to\Ce$,
we have $\D^{A}(f\cdot\phi)=df\cli \phi+f\cdot\D^{A}\phi$,
it follows that the symbol of $\D^{A}$ 
is the Clifford multiplication. Thus $\D^{A}$ is elliptic. 
It has the following \defemph{strong unique continuation property}
(see \cite{arons}, or \cite[ch.~8]{Dirac.boundary}):

\begin{theorem}%
\label{thm-unique.cont}
If $\phi\in\GW$ has a zero of infinite order, and $\D^{A}\phi= 0$ 
for some $A\in\Con(L)$, then $\phi= 0$. 
Consequently, if $\phi\in\GW$ is zero on an open set, and $\D^{A}\phi= 0$ for 
some $A$, then $\phi= 0$. 
In particular, if $\D^{A}\phi=0$, then either $\phi= 0$, or $\supp\phi=M$.
\end{theorem}%

\begin{remark}%
\label{rk-zeros1}
This does \emph{not} mean that the zeros are generic (in this case, isolated). 
In general, the zero set is \emph{countably $2$-rectifiable}
(countable union of Lipschitz images of pieces of a plane), 
so in particular its Hausdorff dimension could be as big as $2$.
See \cite{baer} for a proof.
\end{remark}
The Dirac operator depends on the connection $A\in\Con(L)$ as follows:
\[ \D^{A+2i\theta}\phi=\D^{A}\phi+i\theta \cli\phi \]

\begin{remark}%
\label{rk-SW}
The \defemph{\SW\ equations} are:
$D^{A}\phi=0$ and $\sigma(\phi)=F_{A}^{+}$,
where $iF_{A}$ is the curvature $2$-form of $A$.
See \cite{donaldsonSW} or \cite{morgan} for an introduction to this 
theory.
\end{remark}%

The gauge group $\Ga(L)$ of the unitary bundle $L$ is the group of its 
automorphisms. It can be 
described as
$\Ga(L)=\bigl\{ g:M\to \Sph{1} \bigr\}$, 
acting on $L$ by obvious multiplication.
Its action on $L$ induces an action on $\Con(L)$ by
$g\cdot A=A+2 gdg^{-1}$.
The gauge group $\Ga(L)$ also acts on $\W^{\pm}$ by the obvious
$(g\cdot\phi)\srest{x}=g(x)\cdot\phi\srest{x}$.

Since $\sigma(\phi)=\sigma(\psi)$
\Iff\ $\phi=g\cdot\psi$ for some $g\in\Ga(L)$,
we conclude that the map $\sigma:\GW\to\Gamma(\Lp)$ factors through the orbit space 
$\GW\,\big/\,\Ga(L)$.

The Dirac operator, seen as a map $\GW\times\Con(L)\to 
\Gamma(\W^{-})$, with $(\phi,A)\maps\D^{A}\phi$,
is $\Ga(L)$-invariant as well:
$\D^{g\cdot A}(g\cdot\phi)=g\cdot\D^{A}\phi$.
Therefore, the space of pairs $(\phi,A)$ with $\D^{A}\phi= 0$ is 
a gauge-invariant subspace of $\GW\times\Con(L)$.

\bigskip %

The starting point for proving the results in Table~\ref{table}
are the following simple Lemmata. 
First, a uniqueness result:

\begin{lemma}%
\label{lemma-unique}
If $\phi\in\GW$ admits a connection $A\in\Con(L)$ \st\ 
$\D^{A}\phi= 0$, then $A$ is unique, unless $\phi= 0$.
\end{lemma}
\begin{proof}
Assume $\D^{A}\phi= 0$. If $\phi$ is not constantly-zero, then, 
by the unique continuation property \ref{thm-unique.cont}, 
we must have $\supp\phi=M$.
Suppose $\D^{A+2i\theta}\phi= 0$. 
Since $\D^{A+2i\theta}\phi=\D^{A}\phi+i\theta\cli\phi$, we have
$\theta\cli\phi= 0$. 
But a Clifford product is zero \Iff\ one 
of the factors is zero. 
Thus, since $\supp\phi=M$, we deduce that $\theta=0$ almost-everywhere. By 
continuity, $\theta\equiv 0$.
\end{proof}%

\noindent
And an existence result:

\begin{lemma}%
\label{lemma-exist}
If $\phi\in\GW$ is nowhere-zero, then there is a connection  
$A\in\Con(L)$ \st\ $\D^{A}\phi= 0$.
\end{lemma}
\begin{proof}
Choose a random connection $A$. For any other connection $A+2i\theta$, we have
$\D^{A+2i\theta}\phi=\D^{A}\phi+i\theta\cli\phi$. 
Since $\phi$ is nowhere-zero, 
and since Clifford multiplication is modelled on quaternionic 
multiplication,
the equation $i\theta\cli\phi=-\D^{A}\phi$ can always be solved for 
a $\theta\in\Gamma(\TM^{*})$, and then $\D^{A+2i\theta}\phi= 0$.
\end{proof}%

\noindent
Therefore, each spinor field $\phi$ admits \emph{at most one} 
connection $A$ so that $\D^{A}\phi= 0$. 

We can now easily prove the third line of Table~\ref{table}:

\vbox{%
\begin{lemma}%
$\alpha=\sigma(\phi)$ 
is a bijective correspondence between:
\begin{center}\begin{tabular}{l|l}
Self-dual $2$-forms $\alpha$
&
Gauge classes of pairs $(\phi,A)$
\cr
\quad nowhere-zero
&
\quad with $\D^{A}\phi= 0$ 
\cr
\quad and with $c_{1}(\alpha)=c_{1}(\omega)$.
&
\quad and $\phi$ nowhere-zero.
\end{tabular}\end{center}
\end{lemma}%
}
\begin{proof}
As we assumed that $H^{2}(M;\Zi)$ has no $2$-torsion, the existence of 
a $\phi$ with $\sigma(\phi)=\alpha$ is equivalent to 
$c_{1}(\alpha)=c_{1}(\omega)$. 
Lemmata \ref{lemma-unique} and \ref{lemma-exist} 
insure that there is a unique $A$ \st\ 
$\D^{A}\phi= 0$. And everything is gauge-invariant.
\end{proof}%

\subsection{Splitting the connection}

In what follows, we obtain formulae that completely 
separate the contributions of 
the \LC\ connection $\nabla$  and of $A$ to the spinorial 
connection $\nabla^{A}$ and to the Dirac operator $\D^{A}$. 
These results hold for general 
\spinc-structures, not necessarily associated with an 
\ac\ structure.

The bundle $\W^{+}$ has a natural \herm\  metric, denoted 
by $\Cinner{\cdot,\cdot}$. The real part of this metric is a real 
metric $\Rinner{\cdot,\cdot}=\re\Cinner{\cdot,\cdot}$. In 
what follows, $\inner{\cdot,\cdot}$ will always denote this 
\emph{real} metric, and 
the symbol $\perp$ will denote orthogonality with respect to this real 
metric (and so: $v\perp iv$). Also, all projections will be taken 
with respect to this real metric.

If we change the connection  on $\Ks=\det_{\Ce}\Wp$ from $A $ 
to $A +2i\theta$, then the induced connection on $\W^{+}$ 
changes from $\nabla^{A}$ to 
$\nabla^{A}+i\theta$, and so $\nabla^{A}_{X}\phi$ changes 
to $\nabla^{A}_{X}\phi+i\theta(X)\cdot\phi$. 
Therefore $\nabla^{A}\phi$ changes only in its component along $i\phi$,
while the rest remains fixed, no matter what connection $A $ we 
choose. And indeed:

\begin{proposition}%
\label{prop-nabla.W}
For every \spinc-structure and  connection $A $ 
we have:
\[ \norm{\phi}^{2}\cdot\proj{i\phi^{\perp}}
  \bigl(\nabla^{A} \phi \bigr)
  =i\bigl(\nabla  \sigma(\phi) \bigr)\cli \phi \]
where $\nabla$ is the connection on $\Lp$ induced by the \LC\ 
connection.
\end{proposition}
\begin{proof}
Remember that 
\[ \sigma(\phi)\cli\phi=-i\rec{2}\norm{\phi}^{2}\cdot\phi \]
Apply $\nabla^{A}$: we get 
$\nabla^{A}\bigl( \sigma(\phi)\cli\phi \bigr)
   =-i\rec{2} \nabla^{A}\bigl( \norm{\phi}^{2}\cdot\phi \bigr)$.
Then use the compatibility of the connections with the Clifford 
multiplication:
\[ \bigl( \nabla \sigma(\phi) \bigr)\cli\phi
    +\sigma(\phi)\cli \bigl( \nabla^{A}\phi \bigr)
   =-i\rec{2} d\norm{\phi}^{2}\cdot\phi
    -i\rec{2}\norm{\phi}^{2}\cdot\nabla^{A}\phi \]
and rearrange terms:
\[ \bigl( \sigma(\phi)+i\rec{2}\norm{\phi}^{2} \bigr)\cli\nabla^{A}\phi
   =-i\rec{2} d\norm{\phi}^{2}\cdot\phi
    -\bigl( \nabla \sigma(\phi) \bigr)\cli\phi \]
Since $w\maps\sigma(\phi)\cli w$ is traceless on $\Wp$, it must act 
on each $w\in\Ce\phi$ by:
$w\maps -i\rec{2}\norm{\phi}^{2}\cdot w$, 
while on each $w\in(\Ce\phi)^{\perp}$ it must act by: 
$w\maps +i\rec{2}\norm{\phi}^{2}\cdot w$. 
Therefore $\sigma(\phi)+i\rec{2}\norm{\phi}^{2}$ acts on $\W^{+}$ by 
killing $\Ce\phi$, and by multiplying $(\Ce\phi)^{\perp}$ 
with $i\norm{\phi}^{2}$.
Therefore on $\W^{+}$ we have 
$\sigma(\phi)+i\rec{2}\norm{\phi}^{2}
  =i\norm{\phi}^{2}\cdot\proj{(\Ce\phi)^{\perp}}$, 
and thus:
\[ i\norm{\phi}^{2}
   \cdot\proj{(\Ce\phi)^{\perp}} \bigl(\nabla^{A}\phi \bigr)
   =-i\rec{2} d\norm{\phi}^{2}\cdot\phi
    -\bigl( \nabla \sigma(\phi) \bigr)\cli\phi \]
On the other hand, 
$d\norm{\phi}^{2}=d\Rinner{\phi,\phi}=2\inner{\nabla^{A}\phi,\,\phi}$, 
and therefore 
$d\norm{\phi}^{2}\cdot\phi=2\Inner{\nabla^{A}\phi,\,\phi}\phi
  =2\norm{\phi}^{2}\Inner{\nabla^{A}\phi,\,\rec{\norm{\phi}}\phi}
    \rec{\norm{\phi}}\phi
  =2\norm{\phi}^{2}\cdot\proj{\phi}\bigl(\nabla^{A}\phi \bigr)$.
Consequently:
\[ i\norm{\phi}^{2}
   \cdot\proj{(\Ce\phi)^{\perp}} \bigl(\nabla^{A}\phi \bigr)
   =-i\norm{\phi}^{2}\cdot\proj{\phi}\bigl(\nabla^{A}\phi \bigr)
    -\bigl( \nabla \sigma(\phi) \bigr)\cli\phi \]
and rearranging terms:
\[ i\norm{\phi}^{2} 
   \bigl( \proj{(\Ce\phi)^{\perp}}+\proj{\phi} \bigr) \nabla^{A}\phi
   =-\bigl( \nabla \sigma(\phi) \bigr)\cli\phi \]
Since
$\proj{(\Ce\phi)^{\perp}}+\proj{\phi}=\proj{i\phi^{\perp}}$,
the above is equivalent to the statement of the proposition.
\end{proof}%

\noindent
It should be clear that the formula above (and those following)
can be used carelessly only \emph{off the zeros} of $\phi$.

\begin{corollary} %
For every \spinc-structure and every unitary connection $A $, we 
have: If $\nabla^{A}_{X}\phi=0$, then $\nabla _{X}\sigma(\phi)=0$.
\end{corollary}
\begin{proof}
We need only   remark that, since
$\del_{X}\norm{\phi}^{2}=\del_{X}\Rinner{\phi,\phi}
  =2\inner{\Nabla{X}\phi,\,\phi}$,
if $\nabla^{A}_{X}\phi=0$, then $\del_{X}\norm{\phi}^{2}=0$, and so the length 
of $\phi$ is constant in the $X$-direction. If the length is constantly 0, 
then $\sigma(\phi)=0$ and the statement is trivial. Otherwise, use 
Proposition~\ref{prop-nabla.W} above.
\end{proof}%

\noindent
Remembering \Lemmaref{lemma-kah}, we get:

\begin{corollary}%
\label{cor-kah}
If a non-trivial spinor field $\phi\in\GW$  
admits a connection $A$ \st\ $\nabla^{A}\phi= 0$,
then $\sigma(\phi)$ 
is a \kah\ form compatible with a metric scalar-multiple of $g$.
\end{corollary}%

\begin{remark}%
\label{rk-4k.1}
This last result was obtained independently in \cite{4koreans}. Note 
also that this is only the \spinc\ version of a similar statement 
holding for genuine spin-structures: 
\emph{If $M$ admits a spin-structure 
with a parallel spinor field, then $M$ is \kah.}
The parallel spinor field  reduces 
the holonomy of $M$ to $SU(2)$, and so $M$ is also Ricci-flat, see 
\cite{hitchin.harmonic.spinors} (or 9.18 in \cite[p.~344]{spingeom}).
Because of the extra $ \Sph{1}$-freedom of \spinc-structures, 
in our case the holonomy reduces only to $U(2)$.
\end{remark}%

The part of $\nabla^{A}\phi$ that was left 
undetermined by \Propref{prop-nabla.W} is its component along $i\phi$. 
But since 
$\norm{\phi}^{2}\proj{i\phi}\nabla^{A}\phi
  =\Rinner{\nabla^{A}\phi,\,i\phi}$, 
we get:

\begin{corollary}%
\label{cor-split.conn}
\quad
$\norm{\phi}^{2} \nabla^{A}_{X} \phi 
  =i\bigl(\Nabla{X}\sigma(\phi) \bigr)\cli\phi 
  +\Rinner{\nabla^{A}_{X}\phi,\,i\phi}\cdot i\phi$
\end{corollary}%

The first term on the right side depends only on the \LC\ connection on 
$M$. The second depends only on the choice of $A$. 
By changing $A$ to $A+2i\theta$, 
the term $\inner{\nabla^{A}_{X}\phi,\,i\phi}$ will change to 
$\inner{\nabla^{A}_{X}\phi,\,i\phi}+ \norm{\phi}^{2}\theta(X)$.
For a fixed nowhere-zero $\phi$, suitable choices of $A$ can make 
$\inner{\nabla^{A}\phi,\,i\phi}$  anything one might want. 
In particular --- zero:

\begin{lemma}%
\label{lemma-minim}
For every \spinc-structure and 
every nowhere-zero 
spinor field $\phi\in\GW$, there is a unique
connection $A$ \st\ the induced connection $\nabla^{A}$ on $\W^{+}$ has
\[ \Rinner{\nabla^{A}\phi,\,i\phi}= 0 \]
and so:
$\norm{\phi}^{2} \nabla^{A}\phi
=i\bigl( \nabla \sigma(\phi) \bigr)\cli\phi$.
\end{lemma}
\begin{proof}
The existence part was argued. Uniqueness: 
If both $A$ and $A+2i\theta$ satisfy the 
above, then $i\theta(X)\cdot\phi= 0$, but since $\supp\phi=M$, 
we must have $\theta= 0$.
\end{proof}
\begin{remark}
In \ref{cor-split.conn}, we keep $\phi$ fixed and vary $A$. Then we can see
that the condition $\inner{\nabla^{A}\phi,\,\phi}\rest{x}=0$ is equivalent to 
the covariant derivative $\nabla^{A}\phi\rest{x}$ being \emph{minimal}.
Thus, we could call a connection $A$ with 
$\inner{\nabla^{A}\phi,\,\phi}\equiv 0$ a \defemph{minimal connection} 
for $\phi$.
\end{remark}

\subsection{Splitting the Dirac operator}

Let $\{e_{k}\}$ be any orthonormal orienting local frame in $\TM$. 
Then: 
\begin{align*}
\norm{\phi}^{2}\D^{A}\phi 
 &=\sum e_{k}\cli\norm{\phi}^{2}\Nabla{e_{k}}\phi\\
 &=\sum e_{k}\cli \bigl(
   i\,\Nabla{e_{k}}\sigma(\phi)\cli\phi 
   +\inner{\nabla^{A}_{e_{k}}\phi,\,i\phi}\cdot i\phi \bigr) \\
 &=i\sum e_{k}\cli\bigl( \Nabla{e_{k}}\sigma(\phi)\cli\phi \bigr)
   +i\sum \inner{\nabla^{A}_{e_{k}}\phi,\,i\phi}e_{k}\cli\phi 
\end{align*}
But Clifford multiplication has the property that
\[ v\cli(\alpha\cli\psi)=(v\wedge\alpha-v\cont\alpha)\cli\phi \]
for every $v\in \TM\equiv \TM^{*}$ and $\alpha\in\Lambda^{*}(M)$, 
where $\cont$ is the interior product 
$(v\cont\alpha)(X)=\alpha(v,X)$
(see \cite[p.~25]{spingeom}). 
Also, since the \LC\ connection 
$\nabla$ is torsion-free, we have that
\[ \sum e_{k}\wedge\Nabla{e_{k}}\alpha=d\alpha  \qquad\qquad
  \sum e_{k}\cont\Nabla{e_{k}}\alpha=-d^{*}\alpha \]
where $d$ is exterior differentiation of forms, and $d^{*}=-*d\,*$ is the formal 
adjoint of $d$ with respect to the metric $g$
(see \cite[p.~123]{spingeom}). 
We obtain: 
\[ \norm{\phi}^{2}\D^{A}\phi
  =i\bigl( (d+d^{*})\sigma(\phi)+\inner{\nabla^{A}\phi,\,i\phi} \bigr)
   \cli\phi \]
where we think of $\inner{\nabla^{A}\phi,\,i\phi}$ as the $1$-form 
$X\maps\inner{\nabla^{A}_{X}\phi,\,i\phi}$. 
 
Starting from the Clifford relation
$-e_{1}e_{2}e_{3}e_{4}\cli\phi=\phi$ for $\phi\in\W^{+}$, 
it is easy to prove  that:

\begin{lemma}%
For every $\beta\in\Lambda^{3}(M)$ and every $\phi\in\W^{+}$, we have
$\beta\cli\phi=-(*\beta)\cli\phi$.
Therefore, for every $\alpha\in\Gamma(\Lp)$ and $\phi\in\GW$,
we have $d\alpha\cli\phi=d^{*}\alpha\cli\phi$
\end{lemma}
Therefore:

\begin{theorem}[The  splitting formula]
\label{thm-split.dirac}
\mbox{}\\
For every \spinc-structure and every connection $A$, we have:    
\[ \norm{\phi}^{2} \D^{A}\phi
  =i\bigl( 2\,d^{*}\sigma(\phi) + 
  \Rinner{\nabla^{A}\phi,\,i\phi} \bigr)
  \cli\phi   \]
\end{theorem}
\begin{remark}
\label{rk-zeros2}
This formula is most useful off the zeros of $\phi$. 
Since the structure of the zeros is rather wild in general 
(see \Rkref{rk-zeros1}), one encounters big difficulties when 
trying to obtain something useable at the zeros.
\end{remark}%

\noindent
Some consequences of \ref{thm-split.dirac} are:

\begin{corollary}%
\label{cor-junk}
If $\phi\in\GW$ is a spinor field, and $A$ is a connection \st\ 
$\inner{\nabla^{A}\phi,\,i\phi}\equiv 0$, then:
$\norm{\phi}^{2} \D^{A}\phi=2i\, d^{*}\sigma(\phi)\cli\phi$.
In particular, if $\D^{A}\phi= 0$, then $\sigma(\phi)$ is closed.
\end{corollary}
\begin{proof}
For the second statement, since
$\D^{A}$ is elliptic, it has the unique continuation property. 
Thus, if $\D^{A}\phi=0$, then either $\phi= 0$ and the conclusion 
is trivial, or $\supp\phi=M$. Then, off the zeros 
of $\phi$, $d^{*}\sigma(\phi)=0$, 
and thus by continuity over all $M$.
And $d^{*}=-*d$ on $\Gamma(\Lp)$.
\end{proof}%

\noindent
Further, remembering \Lemmaref{lemma-symp}:

\begin{corollary}%
\label{cor-sympl}
If a nowhere-zero spinor field $\phi$ 
has a  connection $A$ \st\ $\D^{A}\phi= 0$
and $\Rinner{\nabla^{A}\phi,\,i\phi}= 0$, 
then $\sigma(\phi)$ is a symplectic form on $M$, 
compatible with a metric conformal to $g$. 
\end{corollary}%

\begin{remark}%
\label{rk-4k.2}
A similar statement was proved independently in \cite{4koreans}.
There it was proved that: 
\emph{If $\D^{A}\phi= 0$ and 
$\Cinner{\nabla^{A}\phi,\,\phi}= 0$, then $\sigma(\phi)$ is 
symplectic.}
Since the latter statement uses the complex inner product, 
it is more restrictive; 
notice that the $\phi$'s satisfying the latter statement must 
have constant length, 
while those satisfying the Corollary above need not.
\end{remark}%

Remember from \ref{rk-id} the spinor field $\I$. We have:

\begin{corollary}%
Choose a \spinc-structure associated with some \ac\ structure $\omega$. 
Then for the distinguished spinor field $\I\in\GW$ 
and the unique connection $A$ with 
$\inner{\nabla^{A}\I,\,i\I}\equiv 0$,
we have $\D^{A}\I=0$ \Iff\ $\omega$ is symplectic.
\end{corollary}
\begin{remark}
Since $\I$ has constant length, the condition 
$\inner{\nabla^{A}\I,\,i\I}=0$ is equivalent to 
$\nabla^{A}\I \perp \Ce\I$. 
Thus the above is essentially the same as C.~Taubes' 
Lemma~1 from \cite{taubes1}.
It is the easiest of the starting steps of the long investigation
from \cite{taubes.book}.
There, by using the (perturbed) \SW\ equations $\D^{A}\phi=0$ and 
$\sigma(\phi)=F_{A}^{+}-F_{A_{0}}^{+}+\rec{4}\,r\omega$, 
where $A_{0}$ has $\inner{\nabla^{A_{0}}\I,\,i\I}\equiv 0$, 
and $r$ is a positive parameter, Taubes shows roughly that,
in the \emph{symplectic} case, 
\SW\ solutions correspond to holomorphic curves in $(M,\omega)$. 
\end{remark}%

\section{The Correspondence}
\label{sec-core}

In this section, we prove the results announced in Table~\ref{table}.
Remember from \ref{war-torsion} that
$H^{2}(M;\Zi)$ is assumed to have no $2$-torsion.
Proving is now just a matter of 
gathering what was spread around in the paper:

\vbox{%
\begin{theorem}%
\label{thm-a}
$\alpha=\sigma(\phi)$ 
is a bijective correspondence between:
\begin{center}\begin{tabular}{l|l}
\kah\ forms $\alpha$
&
Gauge classes of pairs $(\phi,A)$
\cr
\quad compatible with a metric scalar-multiple of $g$
&
\quad with $\nabla^{A}\phi= 0$, 
\cr
\quad and with $c_{1}(\alpha)=c_{1}(\omega)$.
&
\quad and $\phi$ not constantly zero.
\end{tabular}\end{center}
\end{theorem}%
}
\begin{proof}
A \kah\ form for a metric scalar multiple of $g$ means a nowhere-zero 
\SD\ $2$-form $\alpha$ with $\nabla\alpha= 0$ (\Lemmaref{lemma-kah}).
By \Thmref{thm-lift}, it can be 
lifted to some $\phi$ \Iff\ $c_{1}(\alpha)=c_{1}(\omega)$. 
From left to right, 
suppose given such a \kah\ form $\alpha$. Lift it to a spinor field 
$\phi$ with $\sigma(\phi)=\alpha$. 
The connection $A$ \st $\inner{\nabla^{A}\phi,\,i\phi}\equiv 0$ 
(\Lemmaref{lemma-minim}) has 
$\norm{\phi}^{2}\nabla^{A}\phi
  =i\bigl( \nabla\sigma(\phi) \bigr)\cli\phi$,
so $\nabla^{A}\phi= 0$. And everything is gauge-invariant.
From right to left, use \Corref{cor-kah}.
\end{proof}
\begin{remark}
The above statement can be widened to account for all \kah\ forms 
compatible with a metric \emph{conformal} to $g$ 
(not just scalar-multiple of $g$). 
In one direction, if $\phi$ is a spinor field and 
$\nabla^{A}$ a connection on $\Wp$ \st\ 
\[ v\cli\nabla^{A}_{w}\phi=w\cli\nabla^{A}_{v}\phi \]
for all $v, w\in \TM$, then $\sigma(\phi)$ is a \kah\ form compatible 
with a metric conformal to $g$. 
But to go in the other direction and achieve a bijection similar to the 
one above, one must admit more connections than just the 
$\nabla^{A}$'s. More precisely, one must first embed $\Wp$ into a 
suitable larger bundle, 
and then accept certain unitary connections $\widetilde{\nabla}$ on it that 
do not necessarily preserve $\Wp$. 
In that setting, there is a bijection between: 
on the one hand, \kah\ forms $\alpha$ compatible with a metric conformal 
to $g$ and with $c_{1}(\alpha)=c_{1}(\omega)$, and, on the other 
hand, gauge classes of pairs $(\phi,\widetilde{\nabla})$ 
(with $\phi\in\GW$ and $\widetilde{\nabla}$ an extended connection)
\st\ 
$v\cli\widetilde{\nabla}^{\phantom{A}}_{w}\phi
=w\cli\widetilde{\nabla}^{\phantom{A}}_{v}\phi$ for 
all $v, w \in \TM$.
The discussion of such a statement belongs to a different set of 
ideas, and thus will not be presented here.
\end{remark}
\vbox{%
\begin{theorem}%
\label{thm-b}
$\alpha=\sigma(\phi)$ 
is a bijective correspondence between:
\begin{center}\begin{tabular}{l|l}
Symplectic forms $\alpha$
&
Gauge classes of pairs $(\phi,A)$, with 
\cr
\quad compatible with a metric conformal to $g$
&
\quad $\D^{A}\phi= 0$ \ and \ $\Rinner{\nabla^{A}\phi,\,i\phi}= 0$
\cr
\quad and with $c_{1}(\alpha)=c_{1}(\omega)$.
&
\quad and $\phi$ nowhere-zero. 
\cr
\end{tabular}\end{center}
\end{theorem}%
}
\begin{proof}
Again, $\alpha$ is nowhere-zero and \SD, so $c_{1}(\alpha)=c_{1}(\omega)$ is 
equivalent to the existence of a $\phi$ with $\sigma(\phi)=\alpha$.
From left to right, given $\alpha$ with $d\alpha=0$, lift it to some 
$\phi$. Choose $A$ \st\ $\inner{\nabla^{A}\phi,\,i\phi}\equiv 0$. 
From \ref{cor-junk}, we have 
$\norm{\phi}^{2}\D^{A}\phi=2\,d^{*}\sigma(\phi)\cli\phi$, and thus 
$\D^{A}\phi= 0$. And everything is gauge-invariant. From right to 
left, use \Corref{cor-sympl}.
\end{proof}

\bigskip %

As an amusing application of some methods from this paper, 
we briefly re-prove below 
a classic property of Betti numbers for \kah\ manifolds.
(We entertain the hope that our better peers might use 
techniques from this paper to also uncover some \emph{new} properties of 
\kah\ or symplectic $4$-manifolds.)

\begin{corollary}%
If $M$ is \kah, then $b_{2}^{+}$ must be odd and $b_{1}$ even.
\end{corollary}
\begin{proof}
Being nowhere-zero, the distinguished spinor field $\I$ 
allows us to identify $\uaR\oplus\Lp$ with $\Wp$ (as real bundles) 
through the Clifford action $\xi\maps\xi\cli\I$. 
Then one can compute
\[ \D^{A}(\xi\cli\I)=\bigl( (d+d^{*})\xi \bigr)\cli\I
  +\sum e_{k}\cli\xi\cli\nabla^{A}_{e_{k}}\I \]
But $\xi=f\oplus\alpha$ has $(d+d^{*})\xi=df+d\alpha+d^{*}\alpha$,
and $d\alpha\cli\I=d^{*}\alpha\cli\I$, and so 
$\bigl( (d+d^{*})\xi \bigr)\cli\I=(df+2\,d^{*}\alpha)\cli\I$.

Assume now that $(M,g,\omega)$ is \kah. 
Then, since $\sigma(\I)=\rec{4}\omega$, 
there is a connection $A_{0}$ \st\ $\nabla^{A_{0}}\I\equiv 0$.
In conclusion:
\[ \D^{A_{0}}\bigl( (f\oplus\alpha)\cli\I \bigr)
  =(df+2\,d^{*}\alpha)\cli\I \]
This identifies $\D^{A_{0}}:\GW\to\Gamma(\W^{-})$ 
with $d\oplus2\,d^{*}:\Gamma(\uaR\oplus\Lp)\to\Gamma(\Lambda^{1})$. 
In particular, their kernels must be isomorphic.
Since $d$ and $d^{*}$ have orthogonal images in 
$\Gamma(\Lambda^{1})$, we conclude that $df+2\,d^{*}\alpha=0$ only 
when $df=0$ and $d^{*}\alpha=0$. Thus, the kernel of 
$d\oplus2\,d^{*}$ has dimension $b_{0}+b_{2}^{+}=1+b_{2}^{+}$. 
On the other hand, $\D^{A_{0}}$ is $\Ce$-linear, 
and thus its kernel must be even-dimensional. Therefore $b_{2}^{+}$ 
must be odd.

And finally, since $M$ admits \ac\ \str s, it has $b_{1}+b_{2}^{+}$ odd, 
and so $b_{1}$ must be even.
\end{proof}%

\subsection*{Acknowledgments:}
We wish to thank Rob Kirby for his infinite patience and stimulating 
lectures and conversations. 
Some intersection techniques from Section~1 are inspired from the proof of 
his generalized adjunction formula.
Thanks are also due to A.~K.~Liu for pointing out an embarrassing 
mistake.



\begin{thebibliography}{BLPR01}

\bibitem[AHS78]{AHS}
M.~F. Atiyah, N.~J. Hitchin, and I.~M. Singer, \emph{Self-duality in
  four-dimensional {R}iemannian geometry}, Proc. Roy. Soc. London Ser. A
  \textbf{362} (1978), no.~1711, 425--461. \MR{80d:53023}

\bibitem[Aro57]{arons}
N.~Aronszajn, \emph{A unique continuation theorem for solutions of elliptic
  partial differential equations or inequalities of second order}, J. Math.
  Pures Appl. (9) \textbf{36} (1957), 235--249. \MR{19,1056c}

\bibitem[B{\"a}r97]{baer}
Christian B{\"a}r, \emph{On nodal sets for {D}irac and {L}aplace operators},
  Comm. Math. Phys. \textbf{188} (1997), no.~3, 709--721. \MR{98g:58179}

\bibitem[BBW93]{Dirac.boundary}
Bernhelm Boo{\ss}-Bavnbek and Krzysztof~P. Wojciechowski, \emph{Elliptic
  boundary problems for {D}irac operators}, Birkh\"auser Boston Inc., Boston,
  MA, 1993. \MR{94h:58168}

\bibitem[BLPR01]{4koreans}
Y.~Byun, Y.~Lee, J.~Park, and J.~S. Ryu, \emph{Constructing the {K}\"ahler and
  the symplectic structures from certain spinors on 4-manifolds}, Proc. Amer.
  Math. Soc. \textbf{129} (2001), no.~4, 1161--1168. \MR{1 707 139}

\bibitem[Don96]{donaldsonSW}
S.~K. Donaldson, \emph{The {S}eiberg-{W}itten equations and $4$-manifold
  topology}, Bull. Amer. Math. Soc. (N.S.) \textbf{33} (1996), no.~1, 45--70.
  \MR{96k:57033}

\bibitem[Hit74]{hitchin.harmonic.spinors}
Nigel Hitchin, \emph{Harmonic spinors}, Advances in Math. \textbf{14} (1974),
  1--55. \MR{50 \#11332}

\bibitem[LM89]{spingeom}
H.~Blaine Lawson, Jr. and Marie-Louise Michelsohn, \emph{Spin geometry},
  Princeton University Press, Princeton, NJ, 1989. \MR{91g:53001}

\bibitem[Mor96]{morgan}
John~W. Morgan, \emph{The {S}eiberg-{W}itten equations and applications to the
  topology of smooth four-manifolds}, Princeton University Press, Princeton,
  NJ, 1996. \MR{97d:57042}

\bibitem[Tau94]{taubes1}
Clifford~Henry Taubes, \emph{The {S}eiberg-{W}itten invariants and symplectic
  forms}, Math. Res. Lett. \textbf{1} (1994), no.~6, 809--822. \MR{95j:57039}

\bibitem[Tau00]{taubes.book}
\bysame, \emph{Seiberg-{W}itten and {G}romov invariants for symplectic
  $4$-manifolds}, International Press, Somerville, MA, 2000, Edited by Richard
  Wentworth. \MR{1 798 809}

\bibitem[Wit94]{witten}
Edward Witten, \emph{Monopoles and four-manifolds}, Math. Res. Lett. \textbf{1}
  (1994), no.~6, 769--796. \MR{96d:57035}

\end{thebibliography}

\providecommand{\bysame}{\leavevmode\hbox to3em{\hrulefill}\thinspace}
\providecommand{\MR}{\relax\ifhmode\unskip\space\fi MR }
\providecommand{\MRhref}[2]{%
  \href{http://www.ams.org/mathscinet-getitem?mr=#1}{#2}
}
\providecommand{\href}[2]{#2}


\end{document}